\documentclass[11pt]{article}
\usepackage{amssymb,enumerate}\usepackage{amsmath,amsthm,amsfonts,amssymb, enumerate}
 \numberwithin{equation}{section}
\newtheorem{theo}{Theorem}[section]
\newtheorem{lemme}[theo]{Lemma}
\newtheorem{propo}[theo]{Proposition}
\newtheorem{conj}{Conjecture}
\newtheorem*{theoA}{Theorem (Batty--Robinson)}
\newtheorem*{theoB}{Theorem (Latrach--Mokhtar-Kharroubi)}

\theoremstyle{definition}
\newtheorem{defi}[theo]{Definition}
\newtheorem{nb}[theo]{Remark}
\newtheorem{exe}[theo]{Example}

\newenvironment{proof2}{\noindent {\bf Proof of Theorem \ref{prin}
:}}{\hfill$\blacksquare$\bigskip}
\newenvironment{preuve}{\noindent {\bf Proof
:}}{\hfill$\blacksquare$\bigskip}
\def\ind#1{\lower5pt\hbox{$\scriptstyle #1$}}

\def \d {\mathrm{d}}

\def \uht {\{U_H(t)\,;\,t \geq 0\}}
\def \essu {\operatornamewithlimits{ess\,sup}}
\def \esin {\operatornamewithlimits{ess\,inf}}

\def \leq {\leqslant}
\def \geq {\geqslant}

\def \O {\Omega}
\def \OV {\Omega \times V}
\def \v {v}
\def \w {w}

\renewcommand{\epsilon}{\varepsilon}
\title%[Semigroup generation properties of streaming operators]
{\bf Semigroup generation properties of streaming operators with
non--contractive boundary conditions.}

\author{{\bf Bertrand Lods} \\
\normalsize
Politecnico di Torino, Dipartimento di Matematica,\\
\normalsize Corso Duca degli Abruzzi, 24,\\
\normalsize 10129 Torino, Italy.\\
\normalsize {\tt lods@calvino.polito.it}}

\begin{document}
\date{}
\bibliographystyle{empty}
\maketitle

{\small \noindent {\bf Abstract.} We present $c_0$--semigroup
generation results for the free streaming operator with abstract
boundary conditions. We recall some known results on the matter
and establish a general theorem (already announced in
\cite{cras}). We motivate our study with a lot of examples and
show that our result applies to the physical cases of Maxwell
boundary conditions in the kinetic theory of gases as well as to
the non---local boundary conditions involved in transport--like
equations from population dynamics. \\

}
%\tableofcontents
%%%%%%%%%%%%%%%%%%%%%%%%%%%%%%%%%%%%%%%%%%%%%%%SECTION 1%%%%%%%%%%%%%%%%%%%%%%%%
%
\section{Introduction}\label{sec1}
In this paper, we investigate the well-posedness of the following
initial--boundary--value problem in $L^p$--spaces $(1 \leq p <
\infty)$
\medskip
\begin{subequations}\label{1}
\begin{equation}\label{1a} \dfrac{\partial \phi}{\partial t}(x,\v,t)+\v
\cdot\nabla_x \phi(x,\v,t)=\mathcal{Q}(\phi)(x,\v,t) \qquad \qquad
(x,\v) \in \Omega \times V, \:t > 0\end{equation}
\begin{equation}\label{1b}\phi(x,\v,0)=\phi_0(x,v) \qquad \qquad (x,\v) \in \Omega
\times V \end{equation}
\begin{equation}\label{1c}
\phi_{|\Gamma_-}(x,\v,t)=H(\phi_{|\Gamma_+})(x,\v,t) \qquad
\qquad (x,\v) \in \Gamma_-, t >0.
\end{equation}\end{subequations}
where $\Omega$ is a smooth open subset of $\mathbb{R}^N$ $(N \geq
1)$, $V$ is the support of a positive Radon measure $\d\mu$ on
$\mathbb{R}^N$ and $\phi_0 \in L^p(\Omega \times V,\d x\d\mu(\v))$
$(1 \leq p < \infty).$ The operator $\mathcal{Q}$ at the
right--hand side of \eqref{1a} is a suitable linear operator on
$L^p(\Omega \times V,\d x\d\mu(\v))$. In \eqref{1c} $\Gamma_-$
(respectively $\Gamma_+$) denotes the incoming (resp. outgoing)
part of the boundary of the phase space $\Omega \times V$
$$\Gamma_{\pm}=\{(x,\v) \;\in \partial \Omega \times V\;;\;\pm \v \cdot
n(x) > 0\}$$ where $n(x)$ stands for the outward normal unit at $x
\in \partial \Omega$. The boundary condition \eqref{1c} expresses
that the incoming flux $\phi_{|\Gamma_-}(\cdot,\cdot,t)$ is
related to the outgoing one $\phi_{|\Gamma_+}$ through a
\textit{linear operator} $H$ that we shall assume to be bounded on
some suitable trace spaces.

The kinetic model \eqref{1} arises in different fields of applied
sciences:
\begin{itemize}
\item {\it Mathematical physics.} In the kinetic theory of gases
or in neutron transport theory, the unknown $\phi(x,\v,t)$
represents the density of particles (neutrons, molecules of gas,
etc) having the position $x \in \O$ and the velocity $\v \in V$ at
time $t \geq 0$. In this case, $\mathcal{Q}(\phi)$ represents the
interaction between particles and the host medium due to
collisions \cite{cerci2, cerci3, williams}. \item {\it
Mathematical biology.} In population dynamics, the variables
$(x,\v)$ do no longer represent the position and velocity but any
other state variables of a given cell populations. In this case
$\phi(x,\v,t)$ is the distribution function of cells having the
state $(x,\v)$ at time $t \geq 0$, $\mathcal{Q}(\phi)$ represents
then the transition from one state to another. We refer to
\cite{webb2} for such transport--like equations in the context of
population dynamics and, more generally, to \cite{bello} for
generalized kinetic models in the applied sciences.
\end{itemize} In the present paper, we will focus our attention on
the influence of the boundary operator $H$ on the well--posedness of
\eqref{1}. We will only consider the so--called
\textit{collisionless form} of \eqref{1}, i. e. we will assume that
$$\mathcal{Q}=0.$$
We adopt here the semigroup framework %(see \cite{vdm} for otherconsiderations)
and the main purpose of this paper is to identify the
\textit{right class of boundary operators} $H$ for which the
free--streaming operator (whose domain includes the boundary
condition \eqref{1c})
$$T_H\phi(x,\v)=-\v \cdot \nabla_x \phi(x,\v) \qquad \qquad (x,\v)
\in \OV$$ generates a $c_0$--semigroup in $L^p(\OV,\d x\d\mu(\v))$
$(1 \leq p < \infty).$ Note that, despite the simple aspect of the
transport equation \eqref{1a}, this question is far from being
trivial whenever $H$ is not a contraction. Actually, it is
well--known that for contractive boundary conditions $\|H\| < 1$,
$T_H$ is a generator of a $c_0$--semigroup of contractions in
$L^p(\OV,\d x\d\mu(\v))$ $(1 \leq p < \infty)$ \cite{proto}. The
case of non--contractive boundary conditions is much more involved
because of the difficulty to control the growth of the flux
$\phi(\cdot,\cdot,t).$ We point out that such boundary conditions
arise naturally in population dynamics. Indeed, in this case the
boundary operator $H$ models the birth--law of the cell population
so that $H$ is \textit{multiplicative}. Typically, for a
proliferating population of cells, during the mitosis, mother cells
undergo fission to give birth to two daughter cells, i. e.
$\|Hu\|=2\|u\|$ $\forall u \geq 0$ \cite{leibo}.\\
\noindent The question of the well--posedness of \eqref{1} has
been already addressed in several recent papers, see for instance
\cite{bomt, boul, stream} and the references therein. We present
in this paper various approaches to answer this question and give
also some new results. More precisely, our aim is to determine
sufficient condition on the boundary operator $H$ for which $T_H$
generates a $c_0$-semigroup in $L^p(\OV,\d x\d\mu(\v))$ $(1 \leq p
< \infty).$ Our main result (Theorem \ref{prin}) (already
announced in \cite{cras}) answer this question in general
$L^p$--spaces with arbitrary $1 \leq p < \infty$ by a {\it
constructive approach}. Actually, our proof consists in deriving,
by an appropriate change of unknown, an evolution problem
equivalent to \eqref{1} and involving {\it contractive boundary
conditions}. Note that the afore--mentioned result on contractive
boundary conditions turns out to be a direct consequence of our
main result. Moreover, known results referring to the so--called
phase space approach (see Section \ref{sec3} for more details)
\cite{bou2,bou3} are also simple corollaries of our main theorem.
We apply our results successfully to the following boundary
conditions arising in practical situations:\begin{itemize}
    \item \textit{Local boundary conditions} of Maxwell--type which are known
to be well--suited to the kinetic theory of gases \cite{cerci2}
and to neutron transport theory \cite{williams}.
    \item \textit{Non--local} boundary conditions as the ones used in
population dynamics. Note that this type of boundary conditions may
be handle thanks to {\it compactness arguments}.
\end{itemize} The
outline of the paper is as follows. In the following section, we
present some of the boundary conditions commonly adopted in the
kinetic theory of gases and in population dynamics. These are the
motivating examples we had in mind to apply our main result. In
section 3, we introduce the functional setting and prove the
classical generation theorem for contractive boundary conditions.
In section \ref{sec3}, we present the so--called {\it phase space
approach}. We begin with the particular case of slab geometry
(section \ref{sec31}) and recall then the general result
\cite{bou2} which identify the class of phase spaces in which
\eqref{1} is well-posed without any assumption on the boundary
operator. After some examples showing that, out of this class of
phase spaces, assumptions on the boundary conditions are needed,
we present our main result Theorem \ref{prin} and show that all
the afore--mentioned results are simple consequences of it.
Finally, in section \ref{sec4} we show that our result applies to
the physical boundary conditions afore--mentioned. In an Appendix,
we propose a brief discussion on the use of Batty and Robinson
Theorem \cite{arendt, batty} in the context of kinetic
theory and we end this paper by some concluding remarks and open problems.\\

\noindent {\bf Acknowledgement:} The author would like to thank
Professor Mokhtar--Kharroubi for his precious help and advices
during the preparation of his Ph.D thesis from which the major
part of the present paper is taken. The author aknowledges support
from the European Community through a Marie Curie Individual
Fellowship.
%%%%%%%%%%%%%%%%%%%%%%%%%%%%%%%%%%%%%%%%%%%%%%%SECTION2%%%%%%%%%%%%%%%%%%%%%%%%
%
%
\section{Examples of boundary conditions}\label{sec3}
%%%%%%%%%%%%%%%%%%%%

We present in this section some examples of boundary conditions
arising in applications. These examples are coming from the
\textit{kinetic theory of gases} or from \textit{population
dynamics}. The main feature of these latter is their
\textit{non--local} character whereas the boundary conditions are
\textit{local} in the kinetic theory of gases.
\subsection{Local boundary conditions}
Let us consider in this section the case of Maxwell--type boundary
conditions which plays a fundamental role in the kinetic theory of
gases (see for instance \cite{cerci3}) and in neutron transport
theory \cite{williams}. For simplicity, we assume throughout this
section that $\d\mu(\cdot)$ is the Lebesgue measure with support $V
\subset \mathbb{R}^N$ $(N \geq 1).$ The natural class of boundary
operators arising in the kinetic theory of gases is the one of
boundary operators {\it local with respect to} $x \in
\partial \Omega$.  Typically, such a boundary operator reads
\begin{equation*}\label{frog}
H(\psi_{|\Gamma_+})(x,\v)=\int_{\{v' \in V;\,v'\cdot n(x) > 0\}}
\psi_{|\Gamma_+}(x,v')\,\text{d}\Pi_{(x,\v)}(v') \qquad (x,\v) \in
\Gamma_-,
\end{equation*}
where, for a. e. $(x,\v) \in \Gamma_-,$
$\mathrm{d}\Pi_{(x,\v)}(\cdot)$ is a non--negative and bounded Radon
measure on $\{\v' \in V\,;\,\v'\cdot n(x) > 0\}.$ Precisely,
d$\Pi_{(x,\v)}(\v')$ is the probability that a particle (molecule of
gas, neutron...) striking the {\it wall} $\partial \Omega$ at the
point $x$ with velocity between $\v'$ and $\v'+\d \v'$ will
re--emerge at (practically) the same point with velocity between
$\v$ and $\v+\d \v$ (see \cite{cerci2, cerci3, williams} for
details). A particularly interesting model is the following.
\begin{exe}\label{maxwel} Let us assume that a fraction $\alpha$ $(0<\alpha<1)$ of
particles undergoes a specular reflection while the remaining
fraction $1-\alpha$ is diffused with the Maxwellian distribution of
the wall $M_{\omega}$:
\begin{equation}\label{maxwellienne}
M_{\omega}(\v)=\frac{1}{(2\pi \theta_0)^{N/2}} \exp
(-\frac{\v^2}{2\theta_0})\qquad \qquad \v \in V,\end{equation}
$\theta_0$ being the temperature of the surface $\partial \Omega$
(which is assumed to be constant). Then
 \begin{multline*}\label{maxw}
\text{d}\Pi_{(x,\v)}(\cdot)=\alpha\,\mathrm{d}\delta(\v'-\v+2(\v\cdot
n(x)n(x))+\\
+(1-\alpha)M_{\omega}(\v)|\v' \cdot n(x)|\mathrm{d}\v', \quad (x,\v)
\in \Gamma_-\end{multline*}
 where d$\delta(\cdot)$ is the usual
Dirac mass centered in $0$. This corresponds to the classical
Maxwell model, commonly adopted in the kinetic theory of gases
\cite{cerci2}. \hfill$\diamond$
\end{exe}

More generally, let us introduce the following definition of regular
reflection boundary conditions, due to A. Palczewski \cite{pal}.

\begin{defi} Let $\mathcal{R} \in \mathcal{L}(L^p_+,L^p_-)$. One say that $\mathcal{R}$ is a regular reflection
boundary operator if there exists a $C^1$--piecewise  mapping
$\mathcal{V} \::\: \Gamma_- \to \mathbb{R}^N$ such that
\begin{enumerate}[i)]
\item For any $(x,\v) \in \Gamma_-$, $(x,\mathcal{V}(x,\v)) \in
\Gamma_+.$
\\
\item $|\mathcal{V}(x,\v)|=|\v|$ for any $(x,\v) \in \Gamma_-$.
\\
\item $|n(x) \cdot \v|=|n(x) \cdot \mathcal{V}| |\det
\dfrac{\partial \mathcal{V}}{\partial \v}(x,\v)|,$ $(x,\v) \in
\Gamma_-.$
\\
\item $\mathcal{V}(x,\lambda\v)=\lambda\mathcal{V}(x,\v)$ for any
$(x,\v) \in \Gamma_-$ and $\lambda > 0.$
\\
\item $\mathcal{R}(\varphi)(x,\v)=\varphi(x,\mathcal{V}(x,\v)) \qquad
\qquad \forall (x,\v) \in \Gamma_-, \:\varphi \in L^p_+.$
\end{enumerate}
\end{defi}

\begin{exe} In practical situations, the most frequently used
regular reflection conditions are
\begin{enumerate}[(a)]
\item the {\it specular reflection boundary
conditions} which corresponds to $$\mathcal{V}(x,\v)=\v-2(\v \cdot
n(x))\,n(x) \qquad \qquad (x,\v) \in \Gamma_-.$$
\item The {\it bounce--back
reflection conditions} for which $\mathcal{V}(x,\v)=-\v$, $(x,\v)
\in \Gamma_-$ and $V$ has to be symmetric with respect to
$0$.\hfill$\diamond$ \end{enumerate}\end{exe}

The main important feature of such boundary operators is that they
are {\it conservative}, i. e., for any regular reflection operator
$\mathcal{R}$:
\begin{equation}\label{conser}
\|\mathcal{R}\varphi\|=\|\varphi\| \qquad \forall \varphi \in
L^p_+.\end{equation}

\begin{defi}\label{ck} We shall say that a boundary operator $H \in
\mathcal{L}(L^p_+,L^p_-)$ is of \textit{Maxwell--type} if
\begin{equation*}
H(\psi_{|\Gamma_+})(x,\v)=K(\psi_{|\Gamma_+})(x,\v)+\mathcal{C}(\psi_{|\Gamma_+})(x,\v)\qquad
\qquad (x,\v) \in \Gamma_-, \end{equation*} with $\mathcal{C} \in
\mathcal{L}(L^p_+,L^p_-)$ given by
\begin{equation*}\label{c}
\mathcal{C}(\psi_{|\Gamma_+})(x,\v)=\alpha(x)\mathcal{R}(\psi_{|\Gamma_+})(x,\v)\end{equation*}
where $\alpha(\cdot) \in L^{\infty}(\partial \Omega)$ is
non--negative, $\mathcal{R}$ is a regular reflection operator, and
\begin{equation*}\label{k} K(\psi_{|\Gamma_+})(x,\v)=\int_{\{\v'
\cdot n(x) > 0\}}h(x,\v,\v')\psi_{|\Gamma_+}(x,\v')|\v' \cdot
n(x)|\d \v', \qquad (x,\v) \in \Gamma_-,\end{equation*} where
$h(\cdot,\cdot,\cdot) \geq 0$ is measurable.\end{defi}
\begin{nb} If $\mathcal{C}=0$, the boundary operator is said to be
diffusive. More generally, the operator $K$ is called the
diffusive--part of $H$.
\end{nb}
\subsection{Non--local boundary conditions}

For transport--like equations arising in population dynamics, the
boundary conditions are no longer assumed to be local with respect
to $x \in \partial \Omega$ (see for instance \cite{leibo, webb,
laje} and the monograph \cite{webb2}) as illustrated by the
following example:

\begin{exe}\label{encore} In \cite{bl}, the author, together with M. Mokhtar--Kharroubi,
studied a model of growing cell population proposed by Lebowitz and
Rubinow \cite{leibo}:
\begin{equation}\label{leib}
\begin{cases}\dfrac{\partial \varphi}{\partial t}(a,\ell,t)
+\dfrac{\partial \varphi}{\partial a}(a,\ell,t)+\mu(a,\ell)\varphi(a,\ell,t)=0\\
\varphi(0,\ell)=\displaystyle \int_{\ell_1}^{\ell_2}k(\ell,\ell')\varphi(\ell',\ell')\d\ell'+c\,\varphi(\ell,\ell)\\
 \varphi(a,\ell,0)=\varphi_0(a,\ell) \in X_p
\end{cases}
\end{equation}
where $$\Omega=\{(a,\ell) \in \mathbb{R}^2\,;\,0<a< \ell,\; \ell_1 <
\ell < \ell_2\}$$ with $\mu(\cdot,\cdot) \in L^{\infty}(\Omega)$.
This is a model of a proliferating cell population with inherited
properties. The variable $\ell$ is the cycle length of cells, that
is the time between cell birth and cell division. It is assumed to
be determined at birth. The variable $a$ represents the age of the
individual cell. At birth, the age is obviously null whereas, at
division, $a=\ell.$ The constant $\ell_1$ (respectively $\ell_2$)
denotes the minimum (resp. maximum) cycle length. The unknown
$\varphi(a,\ell,t)$ denotes the density of the cell population with
age $a$ and cycle length $\ell$ at time $t \geq 0.$ The function
$\mu(\cdot,\cdot)$ is the rate of cell mortality which is assumed to
be bounded and non--negative. The boundary condition describes the
birth--law (i.e. the transition from mother cycle length to daughter
cycle length). For this model, the velocity space $V$ reduces to the
singleton
$$V=\{(1,0)\},$$
endowed with the Dirac mass centered in $(1,0)$. One has
$X_p=L^p(\Omega,\d a\d\ell)$ $(1 \leq p <\infty)$ and
$\Gamma_-=\{(0,\ell)\,;\ell_1 < \ell < \ell_2\}$ and
$\Gamma_+=\{(\ell,\ell)\,;\,\ell_1 < \ell < \ell_2\}.$ Let us
consider the biological case
$$\ell_1=0.$$ The
free--streaming operator $T_H$ is given then by
\begin{equation*}T_H\varphi(a,\ell):=- \dfrac{\partial
\varphi}{\partial a}(a,\ell),
\end{equation*}
with its usual domain and, in Eq. \eqref{leib}, the boundary
operator $H \in \mathcal{L}(L^p((0,\ell_2)\,,\d\ell))$ $(1 \leq
p<\infty)$ is {\it non--local} with respect to $x=(a,l) \in \Omega$:
$$H(\psi_{|\Gamma_+})(\ell)=
\int_{0}^{\ell_2}k(\ell,\ell')\psi_{|\Gamma_+}(\ell')\d\ell'+c\,\psi_{|\Gamma_+}(\ell)
\qquad 0 < \ell < \ell_2.$$ \hfill $\diamond$
\end{exe}

As suggested by the above example, we can introduce
\textit{non--local Maxwell--type boundary operators}.
\begin{defi}\label{nonloc}
Let $H \in \mathcal{L}(L^p_+,L^p_-)$, $H$ is said to be a
\textit{non--local Maxwell--type boundary operators} if $H$ writes
$$H=K+\mathcal{C},$$ where $\mathcal{C} \in
\mathcal{L}(L^p_+,L^p_-)$ is a contractive boundary operator and $K
\in \mathcal{L}(L^p_+,L^p_-)$ is a {\it non--local} integral
operator.\end{defi}
\section{Setting of the problem and the classical case of contractive
boundary conditions}\label{sec2}

Let us first introduce the functional setting we shall use in the
sequel. Let
$$X_p=L^p(\Omega \times V,\d x\d \mu(\v)) \qquad \qquad \qquad1 \leq p <
\infty,$$ where $\Omega$ is a smooth interior (respectively
exterior) domain of $\mathbb{R}^N$ $(N \geq 1)$, i.e., $\Omega$ is
bounded (resp. $\mathbb{R}^N \setminus \Omega$ is bounded). The
boundary of the phase space $\partial\O \times V$ splits as
$$\partial \O \times V=\Gamma_- \cup \Gamma_+ \cup \Gamma_0$$
where $\Gamma_{\pm}=\{(x,\v) \;\in \partial \Omega \times
V\;;\;\pm \v \cdot n(x) > 0\}$ and $\Gamma_0=\{(x,\v) \;\in
\partial \Omega \times V\;;\;\pm \v \cdot n(x) = 0\}$. We will
assume throughout this paper that $\Gamma_0$ is of zero measure
with respect to $\d\gamma(\cdot)\d\mu(\cdot)$, $d\gamma(\cdot)$
being the Lebesgue measure on $\partial \Omega$. We define the
partial Sobolev space
$$W_p=\{\psi \in X_p\,;\,\v \cdot \nabla_x \psi \in X_p\}.$$
Suitable $L^p$--spaces for the traces on $\Gamma_{\pm}$ are
defined as
$$L^p_{\pm}=L^p(\Gamma_{\pm};|\v \cdot n(x)|\d \gamma(x)\d\mu(\v)).$$
For any $\psi \in W_p$, one can define the traces
$\psi_{|\Gamma_{\pm}}$ on $\Gamma_{\pm}$, however these traces do
not belong to $L^p_{\pm}$ but to a certain weighted space
\cite{ces1,ces2}. Let us define
$$\widetilde{W}_p=\{\psi \in W_p\,;\,\psi_{|\Gamma_{\pm}} \in
L^p_{\pm}\}.$$ Let $H$ be a bounded linear operator from $L^p_+$
to $L^p_-$
$$H \in {\mathcal L}(L^p_+,L^p_-) \qquad \qquad 1 \leq p < \infty.$$
The free--streaming operator associated with the boundary
condition $H$ is
\begin{equation*}
\begin{cases}
T_H\::\: D(T_H) \subset &X_p \to X_p\\
&\varphi \mapsto T_H\varphi(x,\v):=-\v \cdot \nabla_x \varphi(x,\v),
\end{cases}
\end{equation*}
with domain $$D(T_H):=\{\psi \in \widetilde{W}_p \;\text{ such
that } \psi_{|\Gamma_-}=H(\psi_{|\Gamma_+})\}.$$ A crucial role
will be played in the sequel by the so--called {\it time of
sojourn} in $\Omega$.
\begin{defi}\label{tempsdevol}
For any $(x,\v) \in \overline{\Omega} \times V,$ define
\begin{equation*}\begin{split}
t(x,\v)&=\sup\{\,t > 0\;;x-s\v \in \Omega,\;\;\forall \,0< s <
t\,\}\\
&=\inf\{\,s > 0\,;\,x-s\v \notin
\Omega\}.\end{split}\end{equation*}
For the sake of convenience, we will set
$$\tau(x,\v):=t(x,\v) \:\:\text { if } (x,\v) \in \partial \Omega \times
V.$$
\end{defi}
From a heuristic point of view, $t(x,\v)$ is the time needed by a
particle having the position $x \in \Omega$ and the velocity $-\v
\in V$ to go out $\Omega$. One notes \cite{voigt} that
$\tau(x,\v)=0 \text{ for any } (x,\v) \in \Gamma_-$ whereas, if
$\v \cdot n(x)
>0$, $\tau(x,\v)
> 0.$ In particular,
$$\{(x,\v) \in \Gamma_+\,;\,\tau(x,\v)=0\}=\{(x,\v) \in \Gamma_+\,;\,\v
\cdot n(x)=0\}.$$ Moreover, for any $(x,\v) \in \overline{\Omega}
\times V$
$$(x-t(x,\v)\v,\v) \in \Gamma_-.$$
Let us now derive the resolvent of $T_H$. For any $\lambda \in
\mathbb{C}$ such that Re $\lambda
> 0$, define
\begin{equation*}
\begin{cases}
M_{\lambda} \::\:&L^p_- \to L^p_+\\
&u \mapsto
M_{\lambda}u(x,\v)=u(x-\tau(x,\v)\v,\v)e^{-\lambda\tau(x,\v)},\:\:\:(x,\v)
\in \Gamma_+\;;
\end{cases}
\end{equation*}

\begin{equation*}
\begin{cases}
B_{\lambda} \::\:&L^p_- \to X_p\\
&u \mapsto B_{\lambda}u(x,\v)=u(x-t(x,\v)\v,\v)e^{-\lambda
t(x,\v)},\:\:\:(x,\v) \in \Omega \;;
\end{cases}
\end{equation*}

\begin{equation*}
\begin{cases}
G_{\lambda} \::\:&X_p \to L^p_+\\
&\varphi \mapsto G_{\lambda}\varphi(x,\v)=\displaystyle
\int_0^{\tau(x,\v)}\varphi(x-s\v,\v)e^{-\lambda s}\d
s,\:\:\:(x,\v) \in \Gamma_+\;;
\end{cases}
\end{equation*}
and
\begin{equation*}
\begin{cases}
C_{\lambda} \::\:&X_p \to X_p\\
&\varphi \mapsto C_{\lambda}\varphi(x,\v)=\displaystyle
\int_0^{t(x,\v)}\varphi(x-t\v,\v)e^{-\lambda t}\d t,\:\:\:(x,\v)
\in \Omega\;.
\end{cases}
\end{equation*}
Thanks to H\"older's inequality, all these operators are bounded
on their respective spaces. More precisely, for any
$\mathrm{Re}\lambda
> 0$
\begin{align*}
\|M_{\lambda}\| &\leq 1, & \|B_{\lambda}\| &\leq (p\,\mathrm{Re}
\lambda
)^{-1/p},\\
\\
\|G_{\lambda}\| &\leq (q\,\mathrm{Re} \lambda )^{-1/q}, &
\|C_{\lambda}\| &\leq (\mathrm{Re
}\lambda)^{-1},\:\:\:\:1/p+1/q=1.
\end{align*}
The resolvent of $T_H$ is given by the following (see for instance
\cite{these}).
\begin{propo}\label{reso}
Let $H \in \mathcal{L}(L^p_+,L^p_-)$ be such that there exists
$\lambda_0$ such that
$$r_{\sigma}(M_{\lambda}H) <1\qquad \qquad \forall \,\mathrm{Re} \lambda >
\lambda_0.$$ Then, for any $\mathrm{Re}\lambda
> \lambda_0$,
\begin{equation}\label{resolvante}
(\lambda-T_H)^{-1}=B_{\lambda}H(I-M_{\lambda}H)^{-1}G_{\lambda}+C_{\lambda}.\end{equation}
\end{propo}

We recall now the well--known generation result concerning
contractive boundary conditions. It can be found in \cite[Theorem
2.2, Chapter XII]{proto} (see also \cite{beals}). We recall here
the proof of this (now classical) result since it will play a
fundamental role in the sequel.
\begin{theo}\label{cont}
Let $H \in \mathcal{L}(L^p_+,L_-^p)$ $(1 \leq p < \infty)$ be such
that $\|H\| < 1.$ Then, $T_H$ generates a contraction
$c_0$--semigroup in $X_p$.
\end{theo}

\begin{preuve}  The proof consists in showing that $T_H$ is
\textit{dissipative}. From Proposition \ref{reso}, one sees first
that $\{\lambda \in \mathbb{C} \,;\, \mathrm{Re}\lambda > 0\}\,
\subset \rho(T_H)$, where $\rho(T_H)$ stands for the resolvent set
of $T_H$ (in particular $T_H$ is closed). Let us now consider the
case $1 < p < \infty$ and let $\psi \in D(T_H)$. Since
$$\v \cdot \nabla_x
(|\psi|^p)(x,\v)=p|\psi|^{p-2}(x,\v)\psi(x,\v)(\v \cdot \nabla_x
\psi(x,\v)),$$ one gets
\begin{equation*}\begin{split}
\langle T_H \psi, |\psi|^{p-2}\psi \rangle &:=\int_{\Omega \times
V} |\psi|^{p-2}(x,\v)\,\psi(x,\v) (-\v \cdot \nabla_x
\psi(x,\v))\d x\d\mu(\v)\\
&= -\frac{1}{p} \int_{\Omega \times V} \v \cdot \nabla_x |\psi|^p
(x,\v) \d x\d\mu(\v).\end{split}
\end{equation*}
Green's identity yields
\begin{equation*}\begin{split}
\langle T_H \psi, |\psi|^{p-2}\psi \rangle &=-\frac{1}{p}
\int_{\partial \Omega \times V} |\psi|^p(x,\v)\, \v \cdot
n(x) \d\gamma(x)\d\mu(\v)\\
&=\frac{1}{p} \int_{\Gamma_-} |\psi_{|\Gamma_-}(x,\v)|^p |\v \cdot
n(x)|
\d \gamma(x)\d\mu(\v)\\
&\phantom{+++
+++}-\frac{1}{p}\int_{\Gamma_+}|\psi_{|\Gamma_+}(x,\v)|^p
|\v \cdot n(x)| \d\gamma(x)\d\mu(\v)\\
&=\frac{1}{p}(\|\psi_{|\Gamma_-}\|^p_{L^p_-}-\|\psi_{|\Gamma_+}\|^p_{L^p_+}).
\end{split}\end{equation*}
Since $H$ is a \textit{contraction} and
$\psi_{|\Gamma_-}=H(\psi_{|\Gamma_+})$, one deduces that
$$\langle T_H \psi, |\psi|^{p-2}\psi \rangle < 0.$$
For $p=1$, one shows in the same way that
$$\langle T_H \psi, \text{sign}\psi \rangle < 0 \qquad \forall \; \psi
\in D(T_H).$$ Now, let $\psi \in D(T_H)$ and $\mathrm{Re} \lambda
>0$ be fixed. Set $\varphi=(\lambda-T_H)\psi$ and denotes
\begin{equation*}
\psi^*=\begin{cases} |\psi|^{p-2}\psi &\text{ if } 1< p< \infty\\
\text{sign} \psi &\text{ if } p=1.\end{cases}\end{equation*} One
has $\mathrm{Re }\lambda \|\psi\|^p=\mathrm{Re} \langle \lambda
\psi, \psi^* \rangle.$ Consequently
\begin{equation*}\begin{split}
\mathrm{Re}\lambda \|\psi\|^p &\leq \mathrm{Re }\langle \lambda
\psi, \psi^* \rangle -
\langle T_H \psi, \psi^* \rangle \\
&= \mathrm{Re} \langle \lambda \psi -T_H \psi, \psi^* \rangle \leq
\|\varphi\| \,\|\psi\|^{p-1}.\end{split}\end{equation*} Therefore,
for any $\mathrm{Re} \lambda
> 0$, $\|\psi\| \leq \|\varphi\| /\mathrm{Re} \lambda$, i. e.
\begin{equation}\label{resodiss}
\|(\lambda-T_H)^{-1}\| \leq \frac{1}{\mathrm{Re} \lambda } \qquad
\qquad (\mathrm{Re }\lambda > 0).\end{equation} The proof follows
then from Lumer--Phillips Theorem.
\end{preuve}

\begin{nb}\label{h1}
Note that, resuming the above arguments, one can easily check that
estimate \eqref{resodiss} remains true if one assumes
$\|H\psi\|=\|\psi\|$ $\forall \psi \in L^p_+$. Indeed, with the
notations of the above proof, $\langle T_H \psi, \psi^* \rangle=0$
for any $\psi \in D(T_H).$ Unfortunately, this is not sufficient
to prove that $T_H$ generates a $c_0$--semigroup in $X_p$ as
illustrated by the following example due to J. Voigt \cite{voigt}.
\end{nb}

\begin{exe} Let us consider a 1D transport model in $L^1$. Define
$$\Omega=]0,1[ \qquad \text{ and } \qquad V=[0,\,+\infty[$$
and assume that $\d\mu(\cdot)$ is the Lebesgue measure on $V$. One
sees that $\Gamma_+=\{1\} \times V$ and $\Gamma_-=\{0\} \times V,$
so that
$L^1_{\pm}=L^1([0,\,+\infty[,\mathbf{v}\mathrm{d}\mathbf{v}).$ Let
us consider the \textit{identical} boundary operator
$$H(\psi(1,\cdot))=\psi(0,\cdot) \qquad \forall \psi \in W_1.$$
%It is easy to show that $T_H$ is \textit{not a closed operator} in
%$X_1$ so that $T_H$ is not a generator of a $c_0$--semigroup in
%$X_1$ (see \cite{voigt} or \cite{these} for details).
Let us prove that $T_H$ is not a closed operator in $X_1$. Let $h
\in L^1([0,\,+\infty[,\d\v)$ be such that
\begin{equation}\label{hinfty}
\int_0^{\infty} |h(\v)|\v\d\v=\infty.
\end{equation}
For any $n \in \mathbb{N}$, denote
\begin{equation*}
\varphi_n(x,\v)=\begin{cases} h(\v) &\text{ if } 0< \v < n\\
0 &\text{ else}.\end{cases}\end{equation*} Clearly, $\varphi_n \in
W_1$ for any $n \in \mathbb{N}$ and, since
$$\int_{0}^n|h(\v)|\v\d\v < \infty \qquad \forall n \in
\mathbb{N},$$ one has $\varphi_{n|\Gamma_{\pm}} \in L^1_{\pm}$ and
$\varphi_n \in D(T_H)$ for any $n \in \mathbb{N}$. Now, one can
easily show that
$$\varphi_n \to \varphi \qquad \text{ and } \qquad T_H\varphi_n \to 0
\qquad (n \to \infty)$$ with $\varphi(x,\v)=h(\v)$ for almost
every $(x,v) \in \Omega \times V$, $\varphi \in X_1.$ Now,
according to \eqref{hinfty}
$$\varphi_{\Gamma_-}=h \notin L^1_-.$$
This proves that $\varphi \notin D(T_H)$ and $T_H$ is not a closed
operator in $X_1.$ \hfill$\diamond$
\end{exe}

\begin{nb}\label{V=1} The above example shows that, for $\|H\|=1$, $T_H$ may not be closed and
consequently may not be the generator of a $c_0$-semigroup in $X_p$.
Nevertheless, under the additional assumption $H \geq   0 $, it is
possible to show, by a monotone convergence argument, that there
exists an extension of $T_H$ that generates a $c_0$--semigroup in
$X_p$
 \cite{beals}, \cite[Theorem 2.3, Chapter
XII]{proto}. For more considerations on non--negative conservative boundary conditions, we refer the reader to \cite{arlotti}). %Let us also
%mention the work \cite{mt2} which deals with the Vlasov equation
%in slab geometry with conservative boundary conditions by
%techniques from approximation theory (Chernoff Theorem).
\end{nb}

\begin{nb}\label{h0}
If $\|H\| < 1$, Theorem \ref{cont} implies that the type
$\omega(T_H)$ of the $c_0$-semigroup generated by $T_H$ is
non--positive. Actually, it is possible to derive finer estimates
of $\omega(T_H)$. We refer for instance to \cite{mkak} in the case
when $0 \notin V$ (see also Remark \ref{estimate2} thereafter in
the case of the slab).\end{nb}
%%%%%%%%%%%%%%%%%%%%%%%%%%%%%%%%%%%%%%%%%%%%%%%SECTION3%%%%%%%%%%%%%%%%%%%%%%%%
%
%
\section{The phase space approach}\label{sec3}
%
%%%%%%%%%%%%%%%%%%%%%%%%%%%%%%%%%%%%%%%%%%%%%%%SECTION3.1%%%%%%%%%%%%%%%%%%%%%%%%
%
\subsection{The particular case of a slab}\label{sec31}
We begin this section by dealing with the study of the free
streaming operator in slab geometry. This particular case has its
own historical importance and received a peculiar interest during
the last decade (see for instance \cite{boul,tot}). Precisely, let
$\Omega \times V=]-a,a[ \times [-1,1]$ $(a>0)$ and
$$X_p=L^p(]-a,a[ \times [-1,1],\d x\d\xi) \qquad \qquad (1 \leq p <
\infty).$$ In this case, the incoming and the outgoing part of $\Omega \times V$ are
\begin{equation}\label{gami}\Gamma_-:=\{-a\} \times [0,1] \cup \{a\} \times
[-1,0] \quad \text{ and } \quad \Gamma_+:=\{-a\} \times [-1,0] \cup \{a\} \times
[0,1].\end{equation}
For any $H \in \mathcal{L}(L^p_+,L^p_-)$, the
free streaming operator is given then by
$$T_Hf(x,\xi)=-\xi \dfrac{\partial f}{\partial x}(x,\xi) \qquad f
\in D(T_H)$$ with
$$D(T_H)=\{\psi \in W_p\,;\,H(\psi_{|\Gamma_+})=\psi_{|\Gamma_-}\}, \quad R(T_H) \subset
X_p.$$ It is possible to prove the
following.
\begin{theo}\label{bande}
For any $H \in \mathcal{L}(L^p_+,L^p_-)$, the free streaming
operator $T_H$ is a generator of a $c_0$-semigroup $\{U_H(t)\;;\;t
\geq 0\}$ in $X_p$ $(1 \leq p < \infty).$ Moreover,
\begin{equation}\label{estimate1}
\|U_H(t)\| \leq \max\{1,\|H\|\} \,\exp(t \,
\max\{\frac{1}{2a}\ln{\|H\|},\,0\}),\qquad \qquad t \geq
0.\end{equation}
\end{theo}
\begin{nb}\label{estimate2} Note that, because of the definition of $\Gamma_{\pm}$ \eqref{gami}, any boundary operator $H \in \mathcal{L}(L^p_+,L^p_-)$ admits a
matrix representation \cite{lathe} which allows to improve the estimate \eqref{estimate1} (see \cite{tot}).\end{nb}

This theorem has been proved independently by several authors. Let
us mention here the seminal works of G. Borgioli and S. Totaro
\cite{bomt} and S. Totaro \cite{tot} who proved the result in the
particular case $p=1$ using a general theorem of Batty and
Robinson \cite{batty} (for more details on the result of Batty and
Robinson, see also the Appendix).
%Note that, at least for $p=1$,
%it is possible to improve the estimate \eqref{estimate1} by the
%following: there exists $c \geq 1$ such that
%\begin{equation}\label{estimate2}
%\|U_H(t)\| \leq c \,\exp(t \, \max\{\beta(H),\,0\}),\qquad \qquad
%t \geq 0\end{equation} with
%\begin{equation*}
%\beta(H)=\begin{cases} &\dfrac{1}{2a}\ln\dfrac{2D}{E} \text{ if }
%D\neq
%0\\
%\\
%&\dfrac{1}{2a}\ln{(\|H_{21}\|+\|H_{12}\|)} \text{ if }
%D=0\end{cases}\end{equation*} where
%$D=\|H_{11}\|\|H_{22}\|-\|H_{12}\|\|H_{21}\|$ and
%$$E=\sqrt{(\|H_{12}\|-\|H_{21}\|)^2+4\|H_{11}\|\|H_{22}\|}-(\|H_{21}\|+\|H_{12}\|)$$
%have the same sign \cite{tot}.
More recently, M. Boulanouar proved Theorem \ref{bande} using a
renormalization process similar to that used in Section
\ref{sec3.2} \cite{boul}.\medskip

The above result calls for comments. Surprisingly, Theorem
\ref{bande} asserts that, {\it whatever} the boundary operator $H$
is, the free--streaming operator $T_H$ generates a $c_0$--semigroup
in $X_p$ $(1 \leq p <\infty)$. Actually, as we will see hereafter,
this result follows from the particular nature of the slab geometry.
The drawback of this result is that it does not give any information
of what may occur in other kind of geometry and leaves in the
darkness the real mathematical difficulty. In fact, Theorem
\ref{bande} is a simple consequence of the more general case studied
in the following section.
\subsection{The general phase space approach}\label{sec3.2}

The following illustrates the fact that the geometry of the phase
space plays a crucial role for the well--posedness of kinetic
equations \cite{vdm, bou2}.
\begin{theo}\label{geom} Let us assume that the phase
space $\Omega \times V$ is such that
\begin{equation}\label{minore}
\tau_0:=\esin_{(x,\v) \in \Gamma_+}\tau(x,\v) > 0.\end{equation}
Then, for any $H \in \mathcal{L}(L^p_+,L^p_-)$, $T_H$ is a
generator of a $c_0$--semigroup $\uht$ in $X_p$ $(1 \leq p <
\infty)$ such that
\begin{equation}\label{typeminore}\|U_H(t)\| \leq \max\{1,\|H\|\}
\exp\left(t \max\{0,\ln \|H\|/\tau_0\}\right) \qquad (t \geq
0).\end{equation}
\end{theo}
\begin{nb} Using the terminology of \cite{bou2}, any phase space
$\Omega \times V$ satisfying \eqref{minore} is said to be
\textit{regular}.\end{nb} Transport equations in slab geometry are
governed by the above Theorem since the phase space $[-a,a] \times
[-1,1]$ is \textit{regular}. Indeed, for any $(x,\xi) \in [-a,a]
\times [-1,1]$
\begin{equation*}
t(x,\xi)=\begin{cases} \inf \{s > 0\,;\,x-\xi s \leq -a\} \text{ if } \xi > 0\\
\inf \{s > 0\,;\,x-\xi s \geq a\} \text{ if } \xi <
0,\end{cases}\end{equation*} i. e.
\begin{equation*}\label{Vtempssej1}
t(x,\xi)=\dfrac{x-\text{sign}(\xi)a}{|\xi|}\qquad \qquad (x,\xi) \in [-a,a] \times [-1,1],\;\xi \neq 0.\end{equation*}
Therefore
\begin{equation}\label{tpssej1}
\tau_0=\inf_{(x,\xi) \in D^o} \tau(x,\xi)=2a > 0\end{equation}
which proves that the phase space is {\it regular}.

\begin{nb}\label{coroslab} One notes that estimate \eqref{estimate1} follows from
\eqref{tpssej1} and \eqref{typeminore}. Consequently, Theorem
\ref{bande} turns out to be a simple consequence of Theorem
\ref{geom}. \end{nb}
Theorem \ref{geom} has been proved by M.
Boulanouar \cite{bou2} and his proof is based upon a suitable
\textit{renormalization argument}. More precisely, it consists in
studying the problem
\begin{equation*}\begin{cases}
\dfrac{\d \varphi}{\d t}(t)=T_H\varphi(t)\\
\varphi(0)=\phi_0 \in X_p,\end{cases}\end{equation*} in a
\textit{weighted space} $L^p_\omega:=L^p(\Omega \times
V,\omega(x,\v)\d x\d\mu(\v))$ $(1\leq p < \infty)$ where
$\omega(\cdot,\cdot)$ is a suitable nonnegative function such that
$\omega_{|\Gamma_+}=\|H\|^p,$ $\omega_{|\Gamma_-}=1$ and, because
of \eqref{minore},
\begin{equation*}\label{theta}
\essu_{(x,\v) \in \Omega \times V} \omega(x,\v) \leq \|H\|^p
.\end{equation*} This last inequality implies that the norms on
$X_p$ and on $L^p_\omega$ are \textit{equivalent}. The end of the
proof is based on Hille--Yosida theorem applied in $L^p_{\omega}$
and consists in resuming the arguments of the proof of Theorem
\ref{cont}.

Note that the proof of Theorem \ref{geom} in \cite{vdm} is carried
out by the \textit{method of characteristics}, using the fact
that, because of \eqref{minore}, the lengths of the
characteristics curves have a positive lower bound. The proof of
\cite{vdm} also
uses the above renormalization argument.\\

Theorem \ref{geom} illustrates the important fact that the time of
sojourn is the quantity to handle for who wants to deal with the
well--posedness of linear kinetic equations. Unfortunately, in
practical situations, this theorem only applies in the case of
slab geometry (see Remark \ref{coroslab} above) and in some
particular cases from population dynamics (such like the Rotenberg
model with maturation velocity bounded from below \cite{bouem}).
Indeed, for a bounded convex domain $\Omega \subset \mathbb{R}^N$
$(N \geq 2)$, if $V$ is such that
$$\{\v/|\v|,\,\v \in V, \v \neq 0\}=\mathbb{S}^{N-1} \text{ (the unit sphere of $\mathbb{R}^N$) }$$
then, one can easily check that
$$\inf\{\tau(x,\v),\,(x,\v) \in \Gamma_+\}=0,$$
i.e. $\Omega \times V$ is a {\it non--regular phase space}.
%Therefore, Theorem \ref{geom} is not satisfactory in practical
%situations and one has to prove a more general generation result
%(see Theorem \ref{prin} and Corollary \ref{corgeom} below).
%%%%%%%%%%%%%%%%%%%%%%%%%%%%%%%%%%%%%%%%%%%%%%%SECTION3.3%%%%%%%%%%%%%%%%%%%%%%%%
%
%\subsection{The successive reflection method}\label{sec33}
%

%%%%%%%%%%%%%%%%%%%%%%%%%%%%%%%%%%%%%%%%%%%%%%%SECTION4%%%%%%%%%%%%%%%%%%%%%%%%
%
%
%\section{A new approach: the influence of the boundary operator}\label{sec4}
%
%%%%%%%%%%%%%%%%%%%%%%%%%%%%%%%%%%%%%%%%%%%%%%%SECTION4.2%%%%%%%%%%%%%%%%%%%%%%%%
%
\section{The influence of the boundary
operator}\label{sec4}

%%%%%%%%%%%%%%%%%THEOREM PRINCIPAL%%%%%%%%%%%%

The results of Section \ref{sec3.2} illustrate the fact that, to
prove the well--posedness of kinetic equations associated to a
non--contractive boundary operator $H$, the main difficulty relies
on the fact that, for a convex domain $\Omega \subset
\mathbb{R}^N$ with $N
>1$, the \textit{time of sojourn} of particles in $\Omega$ may be
{\it arbitrary small}. Recall that Theorem \ref{geom} asserts
that, for a regular phase space (for which this time of sojourn is
bounded away from zero), no assumption on the boundary operator is
needed. This is no more the case in full generality as it is
illustrated by the following example:

\begin{exe}[Bounce--back reflections]\label{contrex1} Let $\Omega$  be a smooth open and convex
subset of  $\mathbb{R}^N$ $(N \geq 1)$ and let $V=\mathbb{R}^N$ be
endowed with the Lebesgue measure.  %One note that $\Omega \times
%V$ is a {\it non--regular} phase space, i. e.
%$\inf\{\tau(x,\v)\,;\,(x,\v) \in \Gamma_+\}=0.$
Let us consider
the boundary operator:
$$H(\psi)(x,\v)=\alpha\,\psi(x,-\v) \qquad \qquad (x,\v) \in \Gamma_-,\:\psi \in L^p_+$$
with $\alpha > 1.$ Clearly $H \in \mathcal{L}(L^p_+,L^p_-)$ %since $(x,-\v) \in \Gamma_+$ for any $(x,\v) \in \Gamma_-$
and $\|H\|=\alpha > 1.$ In \cite{lods}, the spectrum of the
associated free--streaming operator $T_H$ is investigated and one
can show that
$$\sigma(T_H) \supseteq \overline{\bigcup_{k \in \mathbb{Z}} R_{\text{ess}}(F_k)}$$
where $R_{\text{ess}}(F_k)$ is the essential range of the measurable mapping:
$$F_k\::\:(x,\v) \in \Omega \times V \mapsto F_k(x,\v)=\dfrac{\ln \alpha-2ik\pi}{t(x,\v)+t(x,-\v)}
\qquad (k \in \mathbb{Z}).$$
Consequently,
$$s(T_H):=\sup\{\mathrm{Re}\lambda\,;\,\lambda \in \sigma(T_H)\}=\essu_{(x,\v) \in \Omega \times V}
\dfrac{\ln \alpha}{t(x,\v)+t(x,-\v)} =+\infty.$$ This proves that
the spectrum of $T_H$ is not confined in any left half--plane. In
particular, $T_H$ is not a generator of a $c_0$--semigroup in
$X_p$ $(1 \leq p < \infty).$\hfill$\diamond$
\end{exe}

\begin{nb} It is possible to exhibit similar examples from neutron
transport models with specular reflection conditions \cite{chen}
and for transport--like equations from population dynamics
\cite{bl} (see also Example \ref{encore} below).\end{nb}

The previous example shows that, for a \textit{non--regular} phase
space, some assumption on the boundary operator is needed to prove
that the associated streaming operator generates a $c_0$--semigroup
in $X_p$. Moreover, Theorem \ref{geom} indicates {\it intuitively}
that $T_H$ will be the generator of a $c_0$--semigroup in $X_p$
provided $H$ "does not take too much into account" the set $\{(x,\v)
\in \Gamma_+,\,\tau(x,\v)=0\}$.
%The following result exploits this intuition and can be found in
%\cite{bou4}.
%\begin{theo}\label{bouzero} Let $H \in \mathcal{L}(L^p_+,L^p_-).$ Assume there exists $\tau_0 > 0$ such that
%\begin{equation}\label{zero}
%H\chi_{\{(x,\v) \in \Gamma_+\,;\,\tau(x,v) <
%\tau_0\}}=0\end{equation} then $T_H$ generates a $c_0$-semigroup
%in $X_p$ $(1 \leq p < \infty).$\end{theo}

%\begin{nb} It is possible to estimate the type of the $c_0$-semigroup
%generated by $T_H$ in terms of $\|H\|$ and $\tau_0$ \cite[Theorem
%3.1]{bou4}. Moreover, a similar result can be obtained for $H$
%such that $\chi_{\{(x,\v) \in \Gamma_-\,;\,\tau(x,-v) <
%\tau_0\}}H=0$, $(\tau_0>0)$ \cite[Theorem 3.1]{bou4}.\end{nb}

%Actually, Theorem \ref{bouzero} may be seen somehow as a corollary
%of Theorem \ref{minore}. Its proof is also based upon a suitable
%renormalization argument. We point out, that unfortunately,
%Theorem \ref{bouzero} is not adapted to the study of practical
%boundary operators (see Section \ref{sec4} thereafter). Assumption
%\eqref{zero} is too restrictive and is be relaxed in the main result of this paper (already announced in \cite{cras}). For any
%$\varepsilon > 0$, denotes $\chi_{\varepsilon}$ the multiplication
%operator in $L^p_+$ by the characteristic function of the set
%$\{(x,\v) \in \Gamma_+\;;\;\tau(x,\v) \leq \varepsilon\}$, i.e.
%$\chi_{\varepsilon} \in \mathcal{L}(L^p_+)$ is given by
Let us make more precise what we mean by this. For any
$\varepsilon > 0$, denotes $\chi_{\varepsilon}$ the multiplication
operator in $L^p_+$ by the characteristic function of the set
$\{(x,\v) \in \Gamma_+\;;\;\tau(x,\v) \leq \varepsilon\}$, i.e.
$\chi_{\varepsilon} \in \mathcal{L}(L^p_+)$ is given by
\begin{equation*}
\chi_{\varepsilon}u(x,v)=\begin{cases}u(x,\v) \text{ if }
\tau(x,\v) \leq
\varepsilon\\
0 \text{ else}, \end{cases}
\end{equation*}
for any $u \in L^p_+.$ Our main result is the following.
\begin{theo}\label{prin}
Let $H \in \mathcal{L}(L^p_+,L^p_-).$ If
\begin{equation}\label{HYP}
\operatornamewithlimits{lim\,sup}_{\varepsilon \to 0}
\|H\chi_{\varepsilon}\|_{\mathcal{L}(L^p_+,L^p_-)} <
1,\end{equation} then $T_H$ generates a $c_0$-semigroup
$\{U_H(t)\;;\;t \geq 0\}$ in $X_p$ $(1 \leq p < \infty).$
Moreover, there exists $C \geq 1$ such that
\begin{equation}\label{type}\|U_H(t)\| \leq C\,\exp (t \max
\{\frac{1}{\varepsilon_0}\ln \|H\|,\,0\}), \qquad \forall\,t \geq
0,\end{equation} where $\varepsilon_0=\sup\{\varepsilon
> 0\,;\,\|H\chi_{\varepsilon}\| < 1\}.$
\end{theo}

\begin{nb}\label{tangent} Roughly speaking, assumption \eqref{prin}
is a smallness assumption of $H$ in the neighborhood of $\{(x,\v)
\in \Gamma_+\;;\,\tau(x,\v)=0\}=\{(x,\v) \in \Gamma_+\,;\,\v \cdot
n(x)=0\}.$ This means that the tangential velocities are weakly
taken into account by $H$ {\it regardless} of its norm.\end{nb}

\begin{nb} A particular version of Theorem \ref{prin} has been first proved in \cite{stream} in the case $p=1$
thanks to Batty--Robinson's theorem. Nevertheless, it appears that
the result of \cite{stream} only apply to regular phase--spaces (see
Appendix for details).\end{nb}

%\begin{nb} Note that Theorem \ref{prin} generalizes a similar
%result by M. Boulanouar \cite{bou3} which asserts that $T_H$
%generates a $c_0$--semigroup in $X_p$ provided there exists
%$\tau_0 > 0$ such that $H$ vanishes on $\{(x,\v) \in
%\Gamma_+\;;\,\tau(x,\v) < \tau_0\}.$ Actually, the result of
%\cite{bou3} is somehow a consequence of Theorem \ref{geom} since
%then, the "support" of $H$ becomes a regular phase space.\end{nb}

%\begin{nb} Obviously Theorem \ref{prin} generalizes the above
%result Theorem \ref{bouzero} (whose publication in \cite{bou4} is
%posterior to the one of \cite{cras}).\end{nb}

\begin{nb} Note that it is possible to show, in the spirit of
\cite[Theorem 4.4]{bl}, that $\{U_H(t)\;;\;t \geq 0\}$ depends
continuously on $H \in \mathcal{L}(L^p_+,L^p_-)$ (see \cite{these}
for details).\end{nb}

Let us explain the strategy we follow to prove this result. This
strategy is inspired by a model from population dynamics (see
Example \ref{encore}) studied together with M. Mokhtar--Kharroubi
\cite{bl}. Our aim is to prove that the following evolution
problem
\begin{equation}\label{evol}\begin{cases}
\dfrac{\partial \psi}{\partial t}(x,\v,t)+\v \cdot \nabla_x
\psi(x,\v,t)=0\\
\psi_{|\Gamma_-}=H(\psi_{|\Gamma_+})\\
\psi(x,\v,0)=\psi_0(x,\v),
\end{cases}\end{equation}
where $\psi_0 \in X_p$ $(1 \leq p < \infty),$ is governed by a
$c_0$--semigroup in $X_p$. We make use of a suitable change of
unknown in the spirit of the one used in \cite{bl} (see also
\cite[Chapter XIII]{proto}). This new unknown satisfies then an
equivalent evolution problem (see below \eqref{evol1}) which,
under assumption \eqref{HYP}, involves a contractive boundary
operator.

Let us introduce some useful definitions. For any $0< q<1$, define
the multiplication operator in $L^p_+$ $(1 \leq p < \infty)$:
$$M_q \::\:u \in L^p_+ \mapsto M_qu(x,\v)=q^{\tau_k(x,\v)} u(x,\v) \in
L^p_+,$$ where $\tau_k(x,\v)=\min \{\tau(x,\v);\,k\},$ $(x,\v) \in
\Gamma_+,$ $k$ being any \textit{fixed} positive real number. Let
$B_q$ be defined by
$$B_q \::\:\varphi \in X_p \mapsto B_q \varphi(x,\v)=q^{t_k(x,\v)}
\varphi(x,\v) \in X_p,$$ with $t_k(x,\v)=\min \{t(x,\v);\,k\},$
$(x,\v) \in \overline{\Omega} \times V.$ Since $M_q \in
\mathcal{L}(L^p_+),$ it is possible to define the absorption
operator associated to $HM_q \in \mathcal{L}(L^p_+,L^p_-)$
\begin{equation*}
\begin{cases}
T_{H_q}\::\: &D(T_{H_q}) \subset X_p \to X_p\\
&\varphi \mapsto T_{H_q}\varphi(x,\v):=-\v \cdot \nabla_x
\varphi(x,\v)-\ln q\,\varphi(x,\v)
\end{cases}
\end{equation*}
where
$$D(T_{H_q})=\{ \psi \in \tilde{W}_p
\;;\;\psi_{|\Gamma_-}=HM_q(\psi_{|\Gamma_+})\}.$$ The unbounded
operators $T_H$ and $T_{H_q}$ are related by the following.
\begin{lemme} For any $0< q < 1$, $B_q^{-1} D(T_H)=D(T_{H_q})$ and
$T_H=B_q T_{H_q} B_q^{-1}.$
\end{lemme}

\begin{preuve} Let $0< q < 1$  be \textit{fixed}. One sees easily that
$B_q$ is a continuous bijection from $X_p$ onto itself. Its
inverse is given by
$$B_q^{-1} \::\:\varphi \in X_p \mapsto
B_q^{-1}\varphi(x,\v)=e^{-t_k(x,\v)\,\ln q} \varphi(x,\v) \in
X_p.$$ Note that $B_q^{-1} \in \mathcal{L}(X_p)$ because $\sup \{
t_k(x,\v)\;;\;(x,\v) \in \Omega \times V\} \leq k.$ Now, let
$\varphi \in D(T_H)$ and $\psi=B_q^{-1}\varphi.$ Let us first show
that $\psi \in W_p$. Indeed, for almost every $(x,\v) \in \Omega
\times V$
\begin{equation*}\begin{split}
\v \cdot \nabla_x \psi(x,\v)&=\lim_{s \to 0}
\frac{\psi(x+s\v,\v)-\psi(x,\v)}{s}\\
&=\lim_{s \to 0} \frac{e^{-t_k(x+s\v,\v) \ln q}
\varphi(x+s\v,\v)-e^{-t_k(x,\v) \ln
q}\varphi(x,\v)}{s}.\end{split}\end{equation*} Since, for a. e.
$(x,\v) \in \Omega \times V$,
$$t(x+s\v,\v)=s+t(x,\v) \qquad \forall \; 0 \leq s < t(x,\v),$$
one gets $t_k(x+s\v,\v)=s+t_k(x,\v)$ for any $0 < s < k-t_k(x,\v)$
and
$$\v \cdot \nabla_x \psi(x,\v)=e^{-t_k(x,\v)\ln q}\; \lim_{s \to 0}
\frac{e^{-s \ln q} \varphi(x+s\v,\v)-\varphi(x,\v)}{s}.$$ Using
that $\varphi \in W_p$ one gets
\begin{equation}\label{napsi}
\v \cdot \nabla_x \psi(x,\v)= e^{-t_k(x,\v)\ln q} (-\ln q\;
\varphi(x,\v) + \v \cdot \nabla_x \varphi(x,\v))
\end{equation}
so that $\psi \in W_p$. Moreover, since $t_k(x,\v)=0$ for any
$(x,\v) \in \Gamma_-,$ it is clear that
\begin{equation*}
\varphi_{|\Gamma_-}=\psi_{|\Gamma_-},\end{equation*} and
$$\psi_{|\Gamma_+}(x,\v)=e^{-\tau_k(x,\v) \ln q}
\varphi_{|\Gamma_+}(x,\v),\:\:\:\:\:\:\:(x,\v) \in \Gamma_+.$$
Thus $\psi_{|\Gamma_{\pm}} \in L^p_{\pm}$ and
$$\psi_{|\Gamma_-}=HM_q(\psi_{|\Gamma_+}).$$
This proves that $\psi \in D(T_{H_q})$ i. e.
$$B_q^{-1} D(T_H) \subset D(T_{H_q}).$$
The converse inclusion is proved similarly. Finally, for $\varphi
\in D(T_H)$, according to \eqref{napsi}
\begin{equation*}\begin{split}
T_{H_q}B_q^{-1}\varphi(x,\v)&=-\v \cdot \nabla_x (e^{-t_k(x,\v)
\ln q}
\varphi(x,\v)) - \ln q\,e^{-t_k(x,\v) \ln q} \varphi(x,\v)\\
&= e^{-t_k(x,\v)\ln q} (\ln q\; \varphi(x,\v) - \v \cdot \nabla_x
\varphi(x,\v)).\end{split}\end{equation*} Consequently
$$B_qT_{H_q}B_q^{-1}\varphi(x,\v)= - \v \cdot \nabla_x \varphi (x,\v) =T_H
\varphi(x,\v)$$ which achieves the proof. \end{preuve}

As a consequence, one has the following.
\begin{propo}\label{equiv}
For any $0 < q < 1$, $T_{H_q}$ generates a $c_0$--semigroup
$\{V_{H_q}(t)\;;$ $\;t \geq 0\}$ in $X_p$ if and only if $T_H$ is
a generator of a $c_0$--semigroup $\{U_{H}(t)\;;\;t \geq 0\}$ in
$X_p$ $(1 \leq p < \infty)$. Moreover,
$$U_H(t)=B_q\,V_{H_q}(t)\,B_q^{-1}\qquad \qquad (t \geq 0).$$
\end{propo}
In other words, Proposition \ref{equiv} indicates that the
following evolution problem
\begin{equation}\label{evol1}\begin{cases}
\dfrac{\partial \varphi}{\partial t}(x,\v,t)+\v \cdot \nabla_x
\varphi(x,\v,t)+\ln q\,\varphi(x,\v,t)=0\\
\varphi_{|\Gamma_-}=HM_q(\varphi_{|\Gamma_+})\\
\varphi(x,\v,0)=e^{- t_k(x,\v) \ln q}\psi_0(x,\v),
\end{cases}\end{equation}
is equivalent to problem \eqref{evol} thanks to the change of
variables $$\varphi(x,\v,t)=e^{- t_k(x,\v) \ln q}\psi(x,\v,t).$$
We are now in position to prove Theorem \ref{prin}.\medskip

\begin{proof2} According to Theorem \ref{cont}, it is enough to prove
the result when $\|H\| \geq 1.$ Define $\mathcal{Q}=\{\,0< q <
1\;;\;\|HM_q\| <1\}$. Proposition \ref{equiv} together with
Theorem \ref{cont} assert that if ${\mathcal Q} \neq \varnothing$
then $T_H$ generates a $c_0$--semigroup $\{U_H(t)\;;\;t \geq 0\}$
such that
\begin{equation}\label{V2.7}U_H(t)=B_qV_{H_q}B_q^{-1}\qquad \forall \,
 t \geq 0,\,q \in \mathcal{Q},\end{equation}
where $\{V_{H_q}(t)\;;\;t \geq 0\}$ is the $c_0$--semigroup in
$X_p$ with generator $T_{H_q}$ $(q \in \mathcal{Q})$.

Thanks to assumption \eqref{HYP}, let us {\it fix} $0<\varepsilon
< k$ so that $\|H\chi_{\varepsilon}\| < 1.$ Then, for any $0<q<1$,
\begin{equation*}\begin{split}
\|HM_q\| &\leq \|H\chi_{\varepsilon}\,M_q\| +
\|H\,(I-\chi_{\varepsilon})\,M_q\|\\
&\leq \|H\chi_{\varepsilon}\|
+\|H\|\,\|(I-\chi_{\varepsilon})\,M_q\|.\end{split}\end{equation*}
Moreover
\begin{equation*}\begin{split}
\|(I-\chi_{\varepsilon})M_q\|&=\sup \{\,e^{\tau_k(x,\v) \ln
q}\;;\,(x,\v) \in \Gamma_+\,\text{ and } \tau_k(x,\v) \geq
\varepsilon\}\\
&\leq e^{\varepsilon \ln q}.\end{split}\end{equation*} Consequently,
$$\|HM_q\| \leq \|H\chi_{\varepsilon}\|+\|H\|e^{\varepsilon \ln{q}}$$
and, if
\begin{equation}\label{2.9}
e^{\varepsilon \ln{q}} <
\frac{1-\|H\chi_{\varepsilon}\|}{\|H\|}\end{equation} then $q \in
\mathcal{Q}$. One has then $\mathcal{Q} \neq \varnothing$ and
$T_H$ is a generator of a $c_0$--semigroup $\{U_H(t)\;;\;t \geq
0\}$ in $X_p$. On the other hand, it is clear that
$$\|V_{H_q}(t)\| \leq e^{-\ln{q}\, t}\qquad \qquad \qquad
\forall\,t \geq 0,\;q \in \mathcal{Q},$$ and one checks that
$$\|B_q\| \leq 1\qquad \qquad \text{ and }
\qquad \qquad \|B_q^{-1}\| \leq e^{-k \ln{q}} \leq e^{-\varepsilon
\ln{q}},\:\:\:\:q \in \mathcal{Q}.$$ Then, \eqref{V2.7} implies
$$\|U_H(t)\| \leq e^{-\varepsilon \ln{q}}\,e^{-\ln{q}\, t}\qquad \qquad
\forall\,t \geq 0,\;q \in \mathcal{Q}.$$ One deduces from
\eqref{2.9} the following estimate
$$\|U_H(t)\| \leq \|H\|e^{\ln{(1-\|H\chi_{\varepsilon}\|)}}
\,e^{\frac{t}{\varepsilon}\ln{\|H\|}} \qquad t\geq 0$$ for any
$\,0 < \varepsilon < k$ such that $\|H\chi_{\varepsilon}\| <1, $
which achieves the proof.\end{proof2}

\begin{nb}\label{nblambda} It has been shown above that, provided $H$ fulfills \eqref{HYP},
$$\operatornamewithlimits{lim\,sup}_{q \to 0}\|HM_q\|<1.$$
Therefore, setting $\lambda=-\ln q$, with the notations of Section \ref{sec2} one gets
$r_{\sigma}(M_{\lambda}H) < 1$ for sufficiently large $\lambda$.
\end{nb}

The results of the previous section are now simple corollaries of
Theorem \ref{prin}. Indeed, let us assume that
\begin{equation*}
\tau_0:=\esin_{(x,\v) \in \Gamma_+}\tau(x,\v) > 0.\end{equation*}
Then, for any bounded operator $H \in \mathcal{L}(L^p_+,L^p_-)$,
one has
\begin{equation}\label{min}
\|H\chi_{\varepsilon}\|=\begin{cases} 0 &\text{ if } 0 <
\varepsilon < \tau_0\\
\|H\| &\text{ if } \varepsilon \geq
\tau_0.\end{cases}\end{equation} Therefore, Theorem \ref{geom} follows directly from
Theorem and assumption \eqref{HYP}
is met by any bounded boundary operator $H$. Note also that the
estimate \eqref{typeminore} follows from \eqref{min} and
\eqref{type}.

%%%%%%%%%%%%%%%%%%%%%%%%%%%%%%%%%%%%%%%%%%%%%%%SECTION5%%%%%%%%%%%%%%%%%%%%%%%%
%
%
\section{Application to Maxwell--type boundary conditions}\label{sec5}
%%%%%%%%%%%%%%%%%%%%

We briefly show in this section how the results of the previous
section apply to the boundary conditions described in Section 2.

We begin by the local boundary conditions introduced in Definition
\ref{ck}. For $p=1$, we have the following.
\begin{propo}
Assume $p=1$ and let $H \in \mathcal{L}(L^1_+,L^1_-)$ be a
Maxwell--type boundary operator given by Definition \ref{ck}. If
\begin{equation*}
\lim_{\varepsilon \to 0}
\,\operatornamewithlimits{ess\,sup}_{\tau(x,\v') \leq \varepsilon}
\int_{\{\v'\cdot n(x) > 0\}}h(x,\v,\v')|\v' \cdot n(x)|\d \v' <1
-\operatornamewithlimits{ess\,sup}_{x \in \partial
\Omega}\,\alpha(x),\end{equation*} then $T_H$ generates a
$c_0$--semigroup in $X_1$. \end{propo}

\begin{preuve} It is easy to check that
\begin{multline*}
\|H\chi_{\varepsilon}\|_{\mathcal{L}(L^1_+,L^1_-)} \leq
\operatornamewithlimits{ess\,sup}_{\tau(x,\v') \leq \varepsilon}
\int_{\{\v'\cdot n(x)
> 0\}}h(x,\v,\v')|\v' \cdot n(x)|\d \v'+\\
+\operatornamewithlimits{ess\,sup}_{x \in \partial
\Omega}\,\alpha(x),\:\:\:\:\:\:\varepsilon > 0.\end{multline*}
Then, Theorem \ref{prin} leads to the conclusion. \end{preuve}

When $1 < p< \infty,$ we have the following.
\begin{propo}
Let $1< p< \infty$. Assume $H$ is a Maxwell--type boundary operator
given by Definition \ref{ck}. Moreover, let us assume that
\begin{multline}\label{diffu1a}
\lim_{\varepsilon \to 0} \operatornamewithlimits{ess\,sup}_{x \in
\partial \Omega} \int_{\{\v \cdot n(x) < 0\}}|\v \cdot n(x)|\d \v \times \\
\times
 \left(\int_{\{\v' \cdot n(x)> 0\} \cap \{\tau(x,\v') \leq
\varepsilon\}}|h(x,\v,\v')|^q |\v' \cdot n(x)|\d
\v'\right)^{\frac{p}{q}} \\
< 1-\operatornamewithlimits{ess\,sup}_{x\in
\partial \Omega}\,\alpha(x) \qquad (1/p+1/q=1).\end{multline} Then
$T_H$ is a generator of a $c_0$--semigroup in $X_p$.
\end{propo}
\begin{preuve}
The proof is a direct application of Theorem \ref{prin} and
follows from straightforward calculations (for the details see
\cite{these}).
\end{preuve}

\begin{nb} It is possible to replace assumption \eqref{diffu1a} by
\begin{multline*}
\lim_{\varepsilon \to 0}
\left(\operatornamewithlimits{ess\,sup}_{(x,\v) \in \Gamma_-}
\int_{\{\v' \cdot n(x)> 0\} \cap \{\tau(x,\v') \leq
\varepsilon\}}h(x,\v,\v')|\v' \cdot n(x)|\d \v'\right)^{\frac{1}{q}} \times  \\
\times \left(\operatornamewithlimits{ess\,sup}_{\tau(x,\v') \leq
\varepsilon} \int_{\{\v \cdot n(x)< 0\}}h(x,\v,\v')|\v \cdot
n(x)|\d \v\right)^{\frac{1}{p}}
<1-\operatornamewithlimits{ess\,sup}_{x \in
\partial \Omega}
\alpha(x)\\
\quad (1/p+1/q=1).\end{multline*} \end{nb}

For practical situations (see Example \ref{maxwel}), it is useful
to state the following.
\begin{propo}\label{propmax} Assume $H=K+\mathcal{C}$ with
$\mathcal{C}$ given by Def. \ref{ck} and
\begin{equation*}
K(\psi_{|\Gamma_+})(x,\v)=\beta(x)\int_{\{\v' \cdot n(x) >
0\}}k(\v,\v')\psi_{|\Gamma_+}(x,\v')|\v' \cdot n(x)| \d
\v',\end{equation*} for any $(x,\v) \in \Gamma_-,$ where
$\beta(\cdot) \in L^{\infty}(\partial \Omega)$ is non--negative.
Moreover, if $\:1 < p <\infty$, assume that
\begin{multline}\label{k(.,.)}
\sup_{x \in \partial \Omega}\int_{\{\v \cdot n(x) < 0\}}|\v \cdot
n(x)|\d \v \left(\int_{\{\v'\cdot n(x) \geq 0\}}|k(\v,\v')|^q|\v'
\cdot n(x)|\d \v'\right)^{p/q}<\infty, \end{multline} where
$1/p+1/q=1$. Then $T_H$ generates a $c_0$--semigroup in $X_p$ $(1
\leq p < \infty)$ provided $\underset{x \in \partial
\Omega}{\essu}\alpha(x) < 1.$
\end{propo}

\begin{preuve}
The proof will consist in showing that the diffusive--part $K$ is
such that
\begin{equation}\label{vanish} \lim_{\varepsilon \to
0}\|K\chi_{\varepsilon}\|=0.\end{equation} We will restrict
ourselves with the case $1 < p < \infty,$ the case $p=1$ being
much simple. For any $\varepsilon > 0,$ define
\begin{multline*}
f_{\varepsilon}(x)=\int_{\{\v \cdot n(x) \leq 0\}}|\v \cdot
n(x)|\d \v \times \\
\times \left(\int_{\{\v'\cdot n(x) \geq 0\} \cap \{\tau(x,\v')
\leq \varepsilon\}}|k(\v,\v')|^q|\v' \cdot n(x)|\d
\v'\right)^{p/q} \qquad (x \in
\partial \Omega).\end{multline*} Clearly, for any $0 \leq
\varepsilon < \varepsilon',$
\begin{equation}\label{crois}
 0 \leq f_{\varepsilon}(x) \leq f_{\varepsilon'}(x) \leq f_{0}(x)
\qquad \qquad (x \in \partial \Omega),\end{equation} where
$$f_0(x)=\int_{\{\v \cdot n(x) \leq 0\}}|\v \cdot n(x)|\d \v
\left(\int_{\{\v'\cdot n(x) \geq 0\}}|k(\v,\v')|^q|\v' \cdot
n(x)|\d \v'\right)^{p/q}.$$ Note that $f_0 \in L^{\infty}(\Omega)$
according to \eqref{k(.,.)}. Moreover, using the continuity of
$n(\cdot)$ and $\tau(\cdot,\cdot)$ (see \cite{voigt}) it is
possible to show \cite[p. 194--195]{these} that
$f_{\varepsilon}(\cdot)$ is continuous on $\partial \Omega$
$(\varepsilon \geq 0).$ Now, for a. e. $(x,\v) \in \Gamma_-$
\begin{equation*}\begin{split}
\lim_{\varepsilon \to 0}&\int_{\{\v' \cdot n(x) \geq 0\} \cap
\{\tau(x,\v') \leq \varepsilon\}}|k(\v,\v')|^q |\v' \cdot n(x)|\d  \v'\\
&=\int_{\{\v' \cdot n(x) \geq 0\} \cap \{\tau(x,\v')
=0\}}|k(\v,\v')|^q |\v'
\cdot n(x)|\d \v'\\
&=\int_{\{\v' \cdot n(x) = 0\}}|k(\v,\v')|^q |\v' \cdot n(x)|\d
\v'=0.\end{split}\end{equation*} Thus, using \eqref{k(.,.)}
together with the dominated convergence theorem,
$$\lim_{\varepsilon \to 0}f_{\varepsilon}(x)=0 \text{ a. e. } x \in
\partial \Omega.$$
Using \eqref{crois} and the continuity of
$f_{\varepsilon}(\cdot)$, Dini's Theorem yields
$$\underset{\varepsilon \to 0}{\lim}\,\underset{ x \in
\partial \Omega}{\sup} f_{\varepsilon }(x)=0.$$
 Now, since
$$\|H\chi_{\varepsilon}\|_{\mathcal{L}(L^p_+,L^p_-)} \leq
\|\beta\|_{\infty}\|f_{\varepsilon}\|^{1/p}_{\infty}$$ one gets
\eqref{vanish}. Finally, since $\|\mathcal{C}\| \leq \essu_{x \in
\partial \Omega} \alpha(x) < 1,$ Theorem \ref{prin} leads to the
conclusion.\end{preuve}

\begin{nb} The main notable fact of Proposition \ref{propmax} is
that generation occurs for arbitrarily large $\beta(\cdot)$. This
comes from the fact that $\beta(\cdot)$ is only space--dependent
and does not care about the tangential velocities (see Remark
\ref{tangent}).
\end{nb}

\begin{exe} Let us consider the Maxwell model described previously. Precisely, assume that, for any $\psi \in L^p_+$,
\begin{multline*}
H(\psi_{|\Gamma_+})(x,\v)=\alpha(x)\,\psi_{|\Gamma_+}(x,\v - 2 (\v
\cdot
n(x))n(x))\\
+(1-\alpha(x))M_{\omega}(\v)\int_{\{\v' \cdot n(x) \geq
0\}}\psi_{|\Gamma_+}(x,\v')|\v' \cdot n(x)|\d \v',\end{multline*}
where $\alpha \in L^{\infty}(\partial \Omega)$ is non--negative
and $M_{\omega}$ is the Maxwellian of the wall given by
\eqref{maxwellienne}. One easily derive from Proposition
\ref{propmax} that, if
$$\sup_{x \in \partial \Omega}\alpha(x) < 1,$$
then $T_H$ is a generator of a $c_0$--semigroup in $X_p$ $(1 < p<
\infty).$ \hfill$\diamond$ \end{exe}

The case of {\it non--local boundary operators} as described in by
Definition \ref{nonloc} is covered by the following result when
$p=1$.
\begin{theo}\label{wcom}
Let $p=1$. Assume that $H=K+\mathcal{C}$ where $\|\mathcal{C}\|<1$
and $K \in \mathcal{L}(L^1_+,\, L^1_-)$ is given by
$$K(\psi)(x,\v)=\int_{\Gamma_+}\kappa(x,\v,y,\v')\psi(y,\v')|\v'\cdot n(y)|\d\gamma(y)d\mu(\v') \qquad (x,\v) \in \Gamma_-$$
where the kernel $\kappa(\cdot,\cdot,\cdot,\cdot) \geq 0$ is
measurable and $\d \gamma(\cdot)$ is the Lebesgue measure on the
surface $\partial \Omega$. If
$$\limsup_{\varepsilon \to 0} \essu_{\{\tau(y,\v') \leq
\varepsilon\}}\int_{\Gamma_-}\kappa(x,\v,y,\v')|\v\cdot
n(x)|\d\gamma(x)\d\mu(\v') < 1 - \|\mathcal{C}\|,$$ then $T_H$
generates a $c_0$--semigroup in $X_p$.
\end{theo}
\begin{preuve} The proof follows from Theorem \ref{prin} and from the fact that
$$\|K\chi_{\varepsilon}\|_{\mathcal{L}(L^1_+,L^1_-)}=\essu_{\{\tau(y,\v') \leq
\varepsilon\}}\int_{\Gamma_-}\kappa(x,\v,y,\v')|\v\cdot
n(x)|\d\gamma(x)\d\mu(v'),$$ since
$\kappa(\cdot,\cdot,\cdot,\cdot)$ is non--negative. \end{preuve}

For $1 < p < \infty$, one has the following result, based on
\textit{compactness arguments}.
\begin{theo}\label{co}
Let $1 < p< \infty$. Assume that $H=K+\mathcal{C}$ where $K
\::L^p_+ \to L^p_-$ is compact and $\|\mathcal{C}\|<1$, then $T_H$
generates a $c_0$--semigroup in $X_p$.
\end{theo}
\begin{preuve} Note that
$$\|H\chi_{\varepsilon}\| \leq \|K\chi_{\varepsilon}\|+
\|\mathcal{C}\|=\|\chi_{\varepsilon}K^\star\|+\|\mathcal{C}\|\:\:\:\forall
\varepsilon > 0$$ where $K^\star \in {\mathcal L}(L^q_-,L^p_+)$
denotes the dual operator of $K$ $ (1/p+1/q=1)$. Since the
truncation operator $\chi_{\varepsilon}$ goes to zero as
$\varepsilon \to 0$ in the strong operator topology (and
consequently uniformly on any compact subset of $L^q_-$) it
follows from the compactness of $K^\star$ that
$$\underset{\varepsilon \to
0}{\lim}\,\|\chi_{\varepsilon}K^\star\|_{{\mathcal
L}(L^q_-,L^q_+)}=0.$$ Hence $\underset{\varepsilon \to
0}{\limsup}\,\|H\chi_{\varepsilon}\| \leq \|\mathcal{C}\| < 1$
which ends the proof thanks to Theorem \ref{prin}. \end{preuve}

{\flushleft \textbf{Example \ref{encore} (revisited).}} Let us go
back to Example \ref{encore}. Let the boundary operator $H \in
\mathcal{L}(L^p((0,\ell_2)\,,\d\ell))$ $(1 \leq p<\infty)$ by
given by
$$H(\psi_{|\Gamma_+})(\ell)=
\int_{0}^{\ell_2}k(\ell,\ell')\psi_{|\Gamma_+}(\ell')\d\ell'+c\,\psi_{|\Gamma_+}(\ell)
\qquad 0 < \ell < \ell_2.$$ If $p=1$, one deduces from Theorem
\ref{wcom} that, provided $$\lim_{\varepsilon \to 0}\,
(\essu_{\ell' \in (0,\varepsilon)}\int_{0}^{\ell_2}
k(\ell,\ell')\,\d\ell) < 1-c$$ then $T_H$ generates a
$c_0$--semigroup in $X_1$ (see also \cite[Corollary 3.2]{bl}). For
$1 < p< \infty,$ it is also possible to prove the well--posedness
of \eqref{leib} thanks to Theorem \ref{co} under some (natural)
assumption on the transition kernel $k(\cdot,\cdot)$ (see
\cite[Corollary 3.1]{bl} for details). \hfill $\diamond$\medskip

\begin{nb} Note
that, if $\ell_1 > 0,$ the phase space $\Omega \times V$ is regular
so that, thanks to Theorem \ref{geom}, $T_H$ generates a
$c_0$--semigroup in $X_p$ for any $H \in
\mathcal{L}(L^p((\ell_1,\ell_2)\,,\d\ell))$ \cite{bl}.\end{nb}

%%%%%%%%%%%%%%%%%%%%%%%%%%%%%%%%%%%%%%%%%%%%%%SECTION6%%%%%%%%%%%%%%%%%%%%%%%%
%
%
%\section{D\'ependance par rapport \`a l'op\'erateur fronti\`ere}\label{sec6}
%
%%
%%

\section{Concluding remarks}\label{conc}

We gave in this paper an overview of $c_0$--semigroup generation
results for free--streaming operators with abstract boundary
conditions. Actually, we emphasize here that, to our mind, the
right approach is the one explained in Section \ref{sec4} which
consists in dealing with the boundary operators rather than with
the phase space. Indeed, for applications, the phase space is
given {\it a priori} and it appears to us that the interesting
question is to determine, for a given phase space, the class of
boundary operators $H$ such that $T_H$ generates a
$c_0$--semigroup in some suitable $L^p$--space. One saw that this
occurs under some suitable smallness assumption on $H$ in the
vicinity of the tangential velocities. The important feature of
such a result (Theorem \ref{prin}) is that no global assumption on
$H$ is needed. Moreover, already known generation results for {\it
regular} phase space turn out to be simple consequence of our main
result. This comes from the fact that, for this kind of geometry,
the set of tangential velocities is empty. We also emphasize the
fact that Theorem \ref{prin} is well--suited to the study of
transport--like equations with practical boundary conditions
arising in the field of mathematical physics (neutron transport
equations, linear kinetic of gases...) or from population
dynamics.

We point out that, by standard  perturbation arguments, the
results of this paper imply the well--posedness (in the semigroup
sense) of the initial-boundary value problem \eqref{1} given in
Introduction with
\begin{equation*}
\mathcal{Q}(f)(x,\v)= \int_V \kappa(x,\v,\w)f(x,\w)\d\mu(\w)-
\sigma(x,\v) f(x,\v)\end{equation*} Precisely, at least for
$\sigma(\cdot,\cdot) \in L^{\infty}(\Omega \times V)$ and for a
measurable kernel $\kappa(\cdot,\cdot,\cdot)$ such that the
operator
$$\mathcal{K}\::\: \psi(x,\v) \in X_p \mapsto \int_V
\kappa(x,\v,\w)f(x,\w,t)\d\mu(\w)\in X_p$$ is {\it bounded}, then
$T_H+\mathcal{Q}$ generates a $c_0$--semigroup in $X_p$ provided
$T_H$ is. It is an open question to know whether such a result is
still valid for unbounded cross--sections $\sigma$ and
$\mathcal{K}$. Such a question is of relevant interest in the
study of the linearized Boltzmann equation (see \cite{cerci3}).
Hopefully, one should generalize the generation result proposed in
\cite{lo4} (dealing with the absorbing case $H=0$) to more general
boundary conditions. Results in this direction are already known
in the peculiar case of slab geometry \cite{chabi1,chabi2} and,
more generally, for regular phase space \cite{vdm}.

We conclude this section with an interesting conjecture. To our
knowledge, all the existing examples of free--streaming operator
$T_H$ that does not generate a $c_0$-semigroup in $X_p$ $(1\leq p
< \infty)$ are such that the spectrum of $T_H$ does not lie in a
left half-space or that $T_H$ is not closed (see Examples 2.5 or
4.1 for instance). Moreover, one saw that the smallness assumption
on $H$ \eqref{HYP} can be seen as an existence assumption of the
resolvent of $T_H$ for large $\lambda$ (see Remark
\ref{nblambda}). This suggests the following conjecture.
\begin{conj} Let $H \in \mathcal{L}(L^p_+,L^p_-)$ $(1 \leq p <
\infty)$ be a bounded boundary operator. Then, $T_H$ generates a
$c_0$-semigroup in $X_p$ if and only if there exists $\lambda_0 \in
\mathbb{R}$ such that $[\lambda_0, +\infty[ \subset
\rho(T_H)$.\end{conj} Actually, the use of Batty--Robinson Theorem
in $L^1$-space (see the following Appendix) supports us in the
belief that the main difficulty to prove that $T_H$ is a generator
is not to find a suitable estimate on the resolvent of $T_H$ but
rather to prove that this resolvent {\it does exist}. Work is in
progress in this direction.
%%%%%%%%%%%%%%%%%%%%%%%%%%%%%%%%%%%%%%%%%%%%%%%SECTION4.1%%%%%%%%%%%%%%%%%%%%%%%%
%
\section*{Appendix: The Batty--Robinson Theorem}\label{sec41}
\renewcommand{\theequation}{A.\arabic{equation}}

In this section, we say a few words about a useful tool used in
kinetic theory to derive generation theorem in $L^1$--space. The
following abstract result is due to J. K. Batty and D. W. Robinson
\cite{batty} (see also \cite{arendt} for a very elegant proof of
this theorem).

Let $X$ be an ordered Banach space whose positive cone is
generating and normal, i.e., $X=X_+-X_-$ and
$X^\star=X^\star_+-X^\star_-$ where $X_{\pm}$ (respectively
$X^\star_{\pm}$) denote the positive and negative cone in $X$
(resp. in $X^\star$).

An operator $A$ on $X$ is said to be \textit{resolvent positive}
if there exists $\omega \in \mathbb{R}$ such that
$]\omega,+\infty[ \subset \rho(A)$ (the resolvent set of $A$) and
$(\lambda-A)^{-1} \geq 0$ for any $\lambda > \omega.$\\
\begin{theoA} Let $A$ be a densely defined resolvent positive operator in
$X$. If there exists $\lambda_0 > s(A)$ and $c>0$ such that
\begin{equation}\label{bat}
\|(\lambda_0-A)^{-1}\,\varphi\| \geq c\|\varphi\| \qquad
\forall\,\varphi \in X^+,\end{equation} then $A$ is a generator of
a (positive) $c_0$--semigroup in $X$.\end{theoA}
Note that the hypothesis \eqref{bat} requires an inverse estimate
with respect to the Hille--Yosida theorem. Note also that, in
practical situations, the Banach space $X$ is a $L^1$--space.

The use of Batty--Robinson's Theorem in kinetic theory is due to
our knowledge to G. Borgioli and S. Totaro \cite{bomt} in order to
prove Theorem \ref{bande} in a $L^1$--setting. More recently, this
result has been used successfully by several authors \cite{mt1,
stream}. In particular, K. Latrach and M. Mokhtar--Kharroubi
\cite{stream} proved a particular version of Theorem \ref{prin}
for $p=1$:
\begin{theoB} Let us assume that $H$ satisfies \eqref{HYP} and the following additional assumptions:
\begin{equation}\label{H>0}
H \geq 0,\end{equation} and
\begin{equation}\label{multi}
\|H\psi\| \geq \|\psi\| \qquad \qquad \forall \psi \in
L_+^1.\end{equation} Then, $T_H$ generates a $c_0$--semigroup in
$L^1(\Omega \times V).$\end{theoB}

Actually, we already saw that according to Remark \ref{nblambda},
there exists $\lambda_0
> 0$ such that
$$r_{\sigma}(M_{\lambda}H) < 1 \qquad \qquad \forall \lambda > \lambda_0.$$
Now, it suffices to appeal to Proposition \ref{reso} together with
\eqref{H>0} which ensure that, for any $\lambda > \lambda_0$,
$(\lambda-T_H)^{-1}$ exists and is {\it nonnegative}. Let us show
how to derive Estimate \eqref{bat}. We follow the strategy of
\cite[Theorem 5.2]{stream}. Let $\lambda>\lambda_0$ and let $\varphi
\in X_1,$ $\varphi \geq 0.$ Set $\psi=(\lambda-T_H)^{-1}\varphi$ the
nonnegative solution of
\begin{equation*}
\lambda \psi(x,\v) + \v \cdot \nabla_x \psi(x,\v)=\varphi(x,\v)
\qquad (x,\v) \in \Omega \times V.\end{equation*} Integrating with
respect to $x$ and $\v$ together with Green's identity leads to
\begin{multline*}
\lambda \|\psi\| + \int_{\Gamma_+}\psi(x,\v)|\v \cdot
n(x)|\d\gamma(x)\d\mu(\v)-\\
-\int_{\Gamma_-}\psi(x,\v)|\v \cdot n(x)|\d\gamma(x)\d\mu(\v)
=\|\varphi\|\end{multline*} which is noting else but
$\lambda\|\psi\| +
\left(\|\psi_{|\Gamma_+}\|-\|H\psi_{|\Gamma_+}\|\right) =
\|\varphi\|.$ Therefore, thanks to \eqref{multi},
$$\|(\lambda-T_H)^{-1}\varphi\| \geq \dfrac{1}{\lambda}\|\varphi\|$$
which gives the estimate \eqref{bat}.

\begin{nb} The above result of \cite{stream} calls for comments.
Actually, it turns out that the assumptions \eqref{HYP} and
\eqref{multi} are compatible only for \textit{regular
phase--space}. Indeed, let us assume that
$\inf\{\tau(x,v)\,;\,(x,v) \in \Gamma_+\}=0$ and define, for any
$\epsilon > 0$, $$\Gamma_{\epsilon}=\{(x,v) \in \Gamma_+\,;\,
\tau(x,v) \leq \epsilon \}$$ and
$$u_{\epsilon}(x,v)=%\dfrac{1}{|\Gamma_{\epsilon}|}
\chi_{\Gamma_{\epsilon}}(x,v) \qquad (x,v) \in \Gamma_+.$$ According
to Assumption \eqref{multi},
$$\|Hu_{\epsilon}\|\geq\|u_{\epsilon}\|$$
and, since $Hu_{\epsilon}=H_{\epsilon}u_{\epsilon}$ (where we used
the notations of Section \ref{sec4}), this shows that
$$\|H_{\epsilon}\| \geq 1$$ and contradicts Assumption \eqref{HYP}. This fact has not been noticed by the
authors of \cite{stream} and suggests that the Batty--Robinson's
Theorem applies in the kinetic theory only to regular
phase--spaces.\end{nb}

%%%%%%%%%%%%%%%%%%%%%%%%%%%%%%
%                            %
%   Bibliographie            %
%                            %
%%%%%%%%%%%%%%%%%%%%%%%%%%%%%%


{\small
\begin{thebibliography}{100}
%
\bibitem{cras}
\textsc{B. Lods}, A generation theorem for kinetic equations with
non--contractive boundary operators. \textit{C. R. Acad. Sci.
Paris.}, Ser. I {\bf 335} 655--660 (2002).
\bibitem{cerci2}
\textsc{C. Cercignani}, \textit{The Boltzmann equation and its
applications}, Springer--Verlag, New York (1988).
\bibitem{cerci3}
\textsc{C. Cercignani, R. Illner and  M. Pulvirenti}, \textit{The
mathematical theory of dilute gases}, Springer--Verlag, New York
(1994).
\bibitem{williams}
\textsc{M. M. R. Williams}, {\it Mathematical Methods in Particle
Transport Theory,} Butterworth, London (1971).
\bibitem{webb2}
\textsc{G. F. Webb},  {\it Theory of nonlinear age--dependent
population dynamics}, Marcel Dekker, New York (1985).
\bibitem{bello}
\textsc{N. Bellomo and M. Pulvirenti Eds.}, {\it Modeling in
applied sciences: A kinetic theory approach}, Birkh\"{a}user,
Boston (2000).
\bibitem{proto}
\textsc{W. Greenberg, C. Van der Mee and V. Protopopescu},
\textit{Boundary Value Problems in Abstract Kinetic theory},
Birkh\"auser Verlag, Basel (1987).
\bibitem{leibo}
\textsc{J. L. Lebowitz and S. I. Rubinow}, A theory for the age
and generation time distribution of a microbial population. {\it
J. Math. Biol.} {\bf 1} 17--36 (1974).
\bibitem{bomt}
\textsc{G. Borgioli and S. Totaro}, 3D--streaming operator with
multiplying boundary conditions: semigroup generation properties.
\textit{Semigroup Forum} {\bf 55} 110--117 (1997).
\bibitem{boul}
 \textsc{M. Boulanouar}, Le transport neutronique avec des conditions
aux limites g\'en\'erales. \textit{C. R. Acad. Sci. Paris.}, Ser.
I {\bf 329} 121--124 (1999).
\bibitem{stream}
\textsc{K. Latrach and M. Mokhtar--Kharroubi}, Spectral analysis
and generation results for streaming operators with multipliying
boundary conditions. \textit{Positivity} {\bf 3} 273--296 (1999).
\bibitem{bou2}
 \textsc{M. Boulanouar}, Op\'erateur d'advection: Existence d'un
semi--groupe (I). \textit{Transp. Theory Stat. Phys.} {\bf 31}
169--176 (2002).
\bibitem{bou3}
 \textsc{M. Boulanouar}, Op\'erateur d'advection: Existence d'un
semi--groupe (II). \textit{Transp. Theory Stat. Phys.} {\bf 32}
185--197 (2003).
\bibitem{arendt}
\textsc{W. Arendt}, Resolvent Positive Operators. \textit{Proc.
London Math. Soc.} {\bf 54} 321--349 (1987).
\bibitem{batty}
\textsc{J. K. Batty and D. W. Robinson}, Positive one parameter
semigroups on ordered spaces. \textit{Acta Appl. Math.} {\bf 1}
221--296 (1984).
 \bibitem{pal}
\textsc{A. Palczewski}, Velocity averaging for boundary value
problems,  In \textit{Nonlinear kinetic theory and mathematical
aspects of hyperbolic systems}, (Edited by V. Boffi, F. Bampi, G.
Toscani), World Scientific. Series Adv. Math. Sci. {\bf Vol. 9}
(1992).
\bibitem{webb}
\textsc{G. F. Webb}, A model of prolifetaring cell population with
inherited cycle length. \textit{J. Math. Biol.} {\bf 23} 269--282
(1986).
\bibitem{laje}
\textsc{K. Latrach and A. Zeghal}, Existence results for a
boundary value problem arising in growing cell populations.
\textit{Math. Models Methods Appl. Sci.} \textbf{ 13} 1--17
(2003).
\bibitem{bl}
\textsc{B. Lods and M. Mokhtar--Kharroubi}, On the theory of a
growing cell population with zero minimum cycle length. \textit{J.
Math. Anal. Appl.} {\bf 266} 70--99 (2001).
\bibitem{ces1}
{\sc M. Cessenat}, Th\'eor\`emes de traces $L_p$ pour les espaces
de fonctions de la neutronique. \textit{C. R. Acad. Sci. Paris.},
Ser I {\bf 299} 831--834 (1984).
\bibitem{ces2}
{\sc M. Cessenat}, Th\'eor\`emes de traces pour les espaces de
fonctions de la neutronique. \textit{C. R. Acad. Sci. Paris.},
Ser. I {\bf 300} 89--92 (1985).
\bibitem{voigt}
\textsc{J. Voigt}, \textit{Functional analytic treatment of the
initial boundary value problem for collisionless gases},
M\"unchen, Habilitationsschrift (1981).
\bibitem{these}
\textsc{B. Lods}, Th\'eorie spectrale des \'equations
cin\'etiques, Th\`ese de doctorat. Universit\'e de
Franche--Comt\'e (2002).
\bibitem{beals}
\textsc{R. Beals and V. Protopopescu}, Abstract time--dependent
transport equations. \textit{J. Math. Anal. Appl.} {\bf 121}
370--405 (1987).
\bibitem{arlotti}
\textsc{L. Arlotti and B. Lods}, work in progress.
\bibitem{mkak}
\textsc{F. Ammar--Khodja and M. Mokhtar--Kharroubi}, On the
exponetial stability of advection semigroups with boundary
operators. \textit{Math. Mod. Meth. Appl. Sci.} {\bf 8} 95--106
(1996).
\bibitem{tot}
\textsc{S. Totaro}, Study of the free streaming operator in slab
geometry in dependence of the boundary conditions. \textit{Math.
Meth. Appl. Sci.} {\bf 20} 717--736 (1997).
\bibitem{lathe}
{\sc K. Latrach}, Th\'eorie spectrales d'\'equations cin\'etiques.
Th\`ese de doctorat. Universit\'e de Franche--Comt\'e (1992).
\bibitem{vdm}
\textsc{C. Van der Mee}, Time dependent kinetic equations with
collision terms relatively bounded with respect to collision
frequency. \textit{Transp. Theory Stat. Phys.} {\bf 30} 63--90
(2001).
\bibitem{bouem}
 \textsc{M. Boulanouar and H. Emamirad}, A transport equation in cell population dynamics. \textit{Differential Integral
   Equations} {\bf 13} 125--144 (2000).
\bibitem{lods}
\textsc{B. Lods}, On the spectrum of tranport operator with
specular and bounce--back reflections conditions. {\it work in
progess}.
\bibitem{chen}
\textsc{Chen Jun and Yang Ming--Zhu}, Linear transport equation
with specular reflection boundary conditions. \textit{Transp.
Theory Stat. Phys.} {\bf 20} 281--306 (1991).
\bibitem{lo4}
\textsc{B. Lods}, On linear kinetic equations involving unbounded
cross--sections. {\it Math. Methods Appl. Sci.} \textbf{27}
1049--1075 (2004).
\bibitem{chabi1}
{\sc M. Chabi and K. Latrach}, On singular mono-energetic
transport equations in slab geometry. \textit{Math. Methods Appl.
Sci.} {\bf 25} 1121--1147 (2002).
\bibitem{chabi2}
{\sc M. Chabi and K. Latrach}, Singular one-dimensional transport
equations on $L\sb p$-spaces. {\it J. Math. Anal. Appl.} {\bf 283}
 319--336 (2003).
\bibitem{mt1}
\textsc{S. Mancini and S. Totaro}, Solutions of the Vlasov
equation in a slab with source terms on the boundaries.
\textit{Riv. Math. Univ. Parma} {\bf 2} 33--47 (1999).









\end{thebibliography}
}
\end{document}